\documentclass[12pt]{article}
\newcommand{\be}{\begin{equation}}
\newcommand{\ee}{\end{equation}}
\newcommand{\bea}{\begin{eqnarray}}
\newcommand{\eea}{\end{eqnarray}}
\newcommand{\ba}{\begin{array}}
\newcommand{\ea}{\end{array}}
\newcommand{\R}{I\!\!R}

\newcommand{\tp}{I\!\!L}

\newcommand{\bc}{\begin{center}}
\newcommand{\ec}{\end{center}}
\newcommand{\ben}{\begin{enumerate}}
\newcommand{\een}{\end{enumerate}}
\newcommand{\bfi}{\begin{figure}}
\newcommand{\efi}{\end{figure}}

\newcommand{\bq}{\begin{quote}}
\newcommand{\eq}{\end{quote}}
\newcommand{\bqu}{\begin{quotation}}
\newcommand{\equ}{\end{quotation}}
\newenvironment{emphit}{\begin{itemize}}{\end{itemize}}
\newcommand{\bemp}{\begin{emphit}}
\newcommand{\eemp}{\end{emphit}}

\newcommand{\bt}{\begin{tabular}}
\newcommand{\et}{\end{tabular}}

\newtheorem{myth}{Theorem}[section]
\newtheorem{mylem}{Lemma}[section]
\newtheorem{mycor}{Corollary}[section]

\newtheorem{mydef}{Definition}[section]
\newtheorem{myrem}{Remark}[section]

\begin{document}
\date{}
\title{Parastrophic invariance of Smarandache quasigroups
}
\author{T\`em\'it\'op\'e Gb\'ol\'ah\`an Ja\'iy\'e\d ol\'a\thanks{On Doctorate Programme at
the University of Abeokuta, Abeokuta, Nigeria.}
\thanks{All correspondence to be addressed to this author}\\
Department of Mathematics,\\
Obafemi Awolowo University, Ile Ife, Nigeria.\\
jaiyeolatemitope@yahoo.com, tjayeola@oauife.edu.ng} \maketitle

\begin{abstract}
Every quasigroup $(L,\cdot )$ belongs to a set of 6 quasigroups,
called parastrophes denoted by $(L,\pi_i)$, $i\in \{1,2,3,4,5,6\}$.
It is shown that $(L,\pi_i )$ is a Smarandache quasigroup with
associative subquasigroup $(S,\pi_i )~\forall~i\in \{1,2,3,4,5,6\}$
if and only if for any of some four $j\in \{1,2,3,4,5,6\}$,
$(S,\pi_j)$ is an isotope of $(S,\pi_i)$ or $(S,\pi_k)$ for one
$k\in \{1,2,3,4,5,6\}$ such that $i\ne j\ne k$. Hence, $(L,\pi_i )$
is a Smarandache quasigroup with associative subquasigroup $(S,\pi_i
)~\forall~i\in \{1,2,3,4,5,6\}$ if and only if any of the six Khalil
conditions is true for any of some four of $(S,\pi_i )$.
\end{abstract}

\section{Introduction}
The study of the Smarandache concept in groupoids was initiated by
W.B. Vasantha Kandasamy in \cite{phd86}. In her book \cite{phd75}
and first paper \cite{phd83} on Smarandache concept in loops, she
defined a Smarandache loop as a loop with at least a subloop which
forms a subgroup under the binary operation of the loop. Here, the
study of Smarandache quasigroups is continued after their
introduction in Muktibodh \cite{muk1} and \cite{muk2}. Let $L$ be a
non-empty set. Define a binary operation ($\cdot $) on $L$ : if
$x\cdot y\in L~\forall ~x, y\in L$, $(L, \cdot )$ is called a
groupoid. If the system of equations ; $a\cdot x=b$ and $y\cdot a=b$
have unique solutions for $x$ and $y$ respectively, then $(L, \cdot
)$ is called a quasigroup. Furthermore, if $\exists$ a $!$ element
$e\in L$ called the identity element such that $\forall ~x\in L$,
$x\cdot e=e\cdot x=x$, $(L, \cdot )$ is called a loop. It can thus
be seen clearly that quasigroups lie in between groupoids and loops.
So, the Smarandache concept needed to be introduced into them and
studied since it has been introduced and studied in groupoids and
loops. Definitely, results of the Smarandache concept in groupoids
will be true in quasigroup that are Smarandache and these together
will be true in Smarandache loops.

It has been noted that every quasigroup $(L,\cdot )$ belongs to a
set of 6 quasigroups, called adjugates by Fisher and Yates
\cite{phd72}, conjugates by Stein \cite{phd73}, \cite{phd76} and
Belousov \cite{phd77} and parastrophes by Sade \cite{phd74}. They
have been studied by Artzy \cite{phd71} and a detailed study on them
can be found in \cite{phd3}, \cite{phd39} and \cite{phd49}. So for a
quasigroup $(L,\cdot )$, its parastrophes are denoted by
$(L,\pi_i)$, $i\in \{1,2,3,4,5,6\}$ hence one can take $(L,\cdot
)=(L,\pi_1)$. For more on quasigroup, loops and their properties,
readers should check \cite{phd3}, \cite{phd41},\cite{phd39},
\cite{phd49}, \cite{phd42} and \cite{phd75}. Let $(G,\oplus)$ and
$(H,\otimes)$ be two distinct quasigroups. The triple $(A,B,C)$ such
that $A,B,C~:~(G,\oplus)\rightarrow (H,\otimes)$ are bijections is
said to be an isotopism if and only if $xA\otimes yB=(x\oplus
y)C~\forall~x,y\in G$. Thus, $H$ is called an isotope of $G$ and
they are said to be isotopic.

In this paper, it will be shown that $(L,\pi_i )$ is a Smarandache
quasigroup with associative subquasigroup $(S,\pi_i )~\forall~i\in
\{1,2,3,4,5,6\}$ if and only if for any of some four $j\in
\{1,2,3,4,5,6\}$, $(S,\pi_j)$ is an isotope of $(S,\pi_i)$ or
$(S,\pi_k)$ for one $k\in \{1,2,3,4,5,6\}$ such that $i\ne j\ne k$.
Hence, it can be concluded that $(L,\pi_i )$ is a Smarandache
quasigroup with associative subquasigroup $(S,\pi_i )~\forall~i\in
\{1,2,3,4,5,6\}$ if and only if any of the six Khalil conditions is
true for any of some four of $(S,\pi_i )$.

\section{Preliminaries}
\begin{mydef}\label{0:0.5}
Let $(L,\cdot )$ be a quasigroup. If there exists at least a
non-trivial subset $S\subset L$ such that $(S,\cdot )$ is an
associative subquasigroup in $L$, then $L$ is called a Smarandache
quasigroup (SQ).
\end{mydef}

\begin{myrem}
Definition~\ref{0:0.5} is equivalent to the definition of
Smarandache quasigroup in \cite{muk1} and \cite{muk2}.
\end{myrem}

\begin{mydef}\label{par:loop}
Let $(L,\theta )$ be a quasigroup. The 5 parastrophes or conjugates
or adjugates of $(L,\theta )$ are quasigroups whose binary
operations
$\theta^*~,~\theta^{-1}~,~{}^{-1}\theta~,~(\theta^{-1})^*~,~({}^{-1}\theta
)^*$ defined on $L$ are given by :
\begin{description}
\item[(a)] $(L,\theta^*)~: ~y\theta ^*x=z\Leftrightarrow x\theta y=z~\forall~x,y,z\in
L$.
\item[(b)] $(L,\theta^{-1})~:~x\theta ^{-1}z=y\Leftrightarrow x\theta y=z~\forall~x,y,z\in
L$.
\item[(c)] $(L,{}^{-1}\theta )~:~z~{}^{-1}\theta y=x\Leftrightarrow x\theta y=z~\forall~x,y,z\in
L.$.
\item[(d)] $(L,(\theta ^{-1})^*)~:~z(\theta ^{-1})^*x=y\Leftrightarrow x\theta y=z~\forall~x,y,z\in
L$.
\item[(e)] $(L,({}^{-1}\theta )^*)~:~y({}^{-1}\theta )^*z=x\Leftrightarrow x\theta y=z~\forall~x,y,z\in
L$.
\end{description}
\end{mydef}

\begin{mydef}\label{rep:loop}
Let $(L,\theta)$ be a loop.
\begin{description}
\item[(a)] $R_x$ and $L_x$ represent the left and right translation maps in
$(L,\theta )~\forall~x\in L$.
\item[(b)] $R_x^*$ and $L_x^*$ represent the left and right translation maps in
$(L,\theta^*)~\forall~x\in L$.
\item[(c)] ${\cal R}_x$ and ${\cal L}_x$ represent the left and right translation maps in
$(L,\theta^{-1})~\forall~x\in L$.
\item[(d)] ${\R}_x$ and ${\tp}$ represent the left and right translation maps in
$(L,{}^{-1}\theta)~\forall~x\in L$.
\item[(e)] ${\cal R}_x^*$ and ${\cal L}_x^*$ represent the left and right translation maps in
$(L,(\theta^{-1})^*)~\forall~x\in L$.
\item[(f)] ${\R}_x^*$ and ${\tp}^*$ represent the left and right translation maps in
$(L,({}^{-1}\theta)^*)~\forall~x\in L$.
\end{description}
\end{mydef}

\begin{myrem}
If $(L,\theta )$ is a loop, $(L,\theta^*)$ is also a loop(and vice
versa) while the other adjugates are quasigroups.
\end{myrem}

\begin{mylem}\label{0:1}
If $(L,\theta )$ is a quasigroup, then
\begin{enumerate}
\item $R_x^*=L_x~,~L_x^*=R_x~,~{\cal L}_x=L_x^{-1}~,~{\R}_x=R_x^{-1}~,~{\cal R}_x^*=L_x^{-1}~,~{\tp}_x^*=R_x^{-1}~\forall~x\in L$.
\item ${\cal L}_x=R_x^{*-1}~,~{\R}_x=L_x^{*-1}~,~{\cal R}_x^*=R_x^{*-1}={\cal L}_x~,~{\tp}_x^*=L_x^{*-1}={\R}_x~\forall~x\in L$.
\end{enumerate}
\end{mylem}
{\bf Proof}\\
The proof of these follows by using Definition~\ref{par:loop} and
Definition~\ref{rep:loop}.
\begin{description}
\item[(1)] $y\theta^*x=z\Leftrightarrow x\theta y=z\Rightarrow y\theta^*x=x\theta y
\Rightarrow yR_x^*=yL_x\Rightarrow R_x^*=L_x$. Also,
$y\theta^*x=x\theta y\Rightarrow
xL_y^*=xR_y\Rightarrow L_y^*=R_y$.\\
$x\theta^{-1}z=y\Leftrightarrow x\theta y=z\Rightarrow x\theta
(x\theta^{-1}z)=z \Rightarrow x\theta z{\cal L}_x=z\Rightarrow
z{\cal L}_xL_x=z\Rightarrow {\cal L}_xL_x=I$. Also,
$x\theta^{-1}(x\theta y)=y\Rightarrow x\theta^{-1}yL_x=y\Rightarrow
yL_x{\cal L}_x=y \Rightarrow L_x{\cal L}_x=I$.
Hence, ${\cal L}_x=L_x^{-1}~\forall~x\in L$.\\
$z({}^{-1}\theta ) y=x\Leftrightarrow x\theta y=z\Rightarrow
(x\theta y)({}^{-1}\theta ) y=x\Rightarrow xR_y({}^{-1}\theta )
y=x\Rightarrow xR_y{\R}_y=x\Rightarrow R_y{\R}_y=I$. Also,
$(z({}^{-1}\theta ) y)\theta y=z\Rightarrow z{\R}_y\theta
y=z\Rightarrow z{\R}_yR_y=z\Rightarrow
{\R}_yR_y=I$. Thence, ${\R}_y=R_y^{-1}~\forall~x\in L$.\\
$z(\theta^{-1})^*x=y\Leftrightarrow x\theta y=z$, so, $x\theta
(z(\theta^{-1})^*x)=z\Rightarrow x\theta z{\cal R}_x^*=z \Rightarrow
z{\cal R}_x^*L_x=z\Rightarrow {\cal R}_x^*L_x=I$. Also, $(x\theta
y)(\theta^{-1})^*x=y\Rightarrow yL_x(\theta ^{-1})^*x=y\Rightarrow
yL_x{\cal R}_x^*=y\Rightarrow L_y{\cal R}_x^*=I$.
Whence, ${\cal R}_x^*=L_x^{-1}$.\\
$y({}^{-1}\theta )^*z=x\Leftrightarrow x\theta y=z$, so,
$y({}^{-1}\theta )^*(x\theta y)=x\Rightarrow
y({}^{-1}\theta)^*xR_y=x\Rightarrow xR_y{\tp}_y^*=x\Rightarrow
R_y{\tp}_y^*=I$. Also, $(y({}^{-1}\theta)^*z)\theta y=z\Rightarrow
z{\tp}_y^*\theta y=z\Rightarrow z{\tp}_y^*R_y=z\Rightarrow
{\tp}_y^*R_y=I$. Thus, ${\tp}_y^*=R_y^{-1}$.
\item[(2)] These ones follow from (1).
\end{description}

\begin{mylem}\label{0:1.1}
Every quasigroup which is a Smarandache quasigroup has at least a
subgroup.
\end{mylem}
{\bf Proof}\\
If a quasigroup $(L,\cdot )$ is a SQ, then there exists a
subquasigroup $S\subset L$ such that $(S,\cdot )$ is associative.
According \cite{phd2}, every quasigroup satisfying the associativity
law has an identity hence it is a group. So, $S$ is a subgroup of
$L$.

\begin{myth}\label{0:1.2}(Khalil Conditions \cite{phd87})
A quasigroup is an isotope of a group if and only if any one of some
six identities are true in the quasigroup.
\end{myth}

\section{Main Results}
\begin{myth}\label{1:1}
$(L,\theta )$ is a Smarandache quasigroup with associative
subquasigroup $(S,\theta )$ if and only if any of the following
equivalent statements is true.
\begin{enumerate}
\item $(S,\theta )$ is isotopic to $(S,(\theta^{-1})^*)$.
\item $(S,\theta^*)$ is isotopic to $(S,\theta^{-1})$.
\item $(S,\theta)$ is isotopic to $(S,({}^{-1}\theta )^*)$.
\item $(S,\theta^*)$ is isotopic to $(S,{}^{-1}\theta )$.
\end{enumerate}
\end{myth}
{\bf Proof}\\
$L$ is a SQ with associative subquasigroup $S$ if and only if
$s_1\theta (s_2\theta s_3)=(s_1\theta s_2)\theta s_3\Leftrightarrow
R_{s_2}R_{s_3}=R_{s_2\theta s_3}\Leftrightarrow L_{s_1\theta
s_2}=L_{s_2}L_{s_1}~\forall~s_1,s_2,s_3\in S$.

The proof of the equivalence of (1) and (2) is as follows.
$L_{s_1\theta s_2}=L_{s_2}L_{s_1}\Leftrightarrow {\cal L}_{s_1\theta
s_2}^{-1}={\cal L}_{s_2}^{-1}{\cal L}_{s_1}^{-1}\Leftrightarrow
{\cal L}_{s_1\theta s_2}={\cal L}_{s_1}{\cal L}_{s_2}\Leftrightarrow
(s_1\theta
s_2)\theta^{-1}s_3=s_2\theta^{-1}(s_1\theta^{-1}s_3)\Leftrightarrow
(s_1\theta s_2){\cal R}_{s_3}=s_2\theta^{-1}s_1{\cal
R}_{s_3}=s_1{\cal R}_{s_3}(\theta^{-1})^*s_2\Leftrightarrow
(s_1\theta s_2){\cal R}_{s_3}=s_1{\cal
R}_{s_3}(\theta^{-1})^*s_2\Leftrightarrow (s_2\theta^* s_1){\cal
R}_{s_3}=s_2\theta^{-1}s_1{\cal R}_{s_3}\Leftrightarrow ({\cal
R}_{s_3}, I, {\cal R}_{s_3})~:~(S,\theta )\to
(S,(\theta^{-1})^*)\Leftrightarrow (I, {\cal R}_{s_3},{\cal
R}_{s_3})~:~(S,\theta^* )\to (S,\theta^{-1})\Leftrightarrow
(S,\theta )$ is isotopic to $(S,(\theta^{-1})^*)\Leftrightarrow
(S,\theta^*)$ is isotopic to $(S,\theta^{-1})$.

The proof of the equivalence of (3) and (4) is as follows.
$R_{s_2}R_{s_3}=R_{s_2\theta s_3}\Leftrightarrow
{\R}_{s_2}^{-1}{\R}_{s_3}^{-1}={\R}_{s_2\theta
s_3}^{-1}\Leftrightarrow {\R}_{s_3}{\R}_{s_2}={\R}_{s_2\theta
s_3}\Leftrightarrow (s_1{}^{-1}\theta s_3){}^{-1}\theta
s_2=s_1{}^{-1}\theta (s_2\theta s_3)\Leftrightarrow (s_2\theta
s_3){\tp}_{s_1}=s_3{\tp}_{s_1}{}^{-1}\theta s_2=s_2({}^{-1}\theta
)^*s_3{\tp}_{s_1}\Leftrightarrow (s_2\theta
s_3){\tp}_{s_1}=s_2({}^{-1}\theta )^*s_3{\tp}_{s_1}\Leftrightarrow
(s_3\theta^* s_2){\tp}_{s_1}=s_3{\tp}_{s_1}{}^{-1}\theta
s_2\Leftrightarrow (I, {\tp}_{s_1}, {\tp}_{s_1})~:~(S,\theta )\to
(S,({}^{-1}\theta)^*)\Leftrightarrow ({\tp}_{s_1},I,
{\tp}_{s_1})~:~(S,\theta^* )\to (S,{}^{-1}\theta )\Leftrightarrow
(S,\theta)$ is isotopic to $(S,({}^{-1}\theta )^*)\Leftrightarrow
(S,\theta^*)$ is isotopic to $(S,{}^{-1}\theta )$.

\begin{myrem}
In the proof of Theorem~\ref{1:1}, it can be observed that the
isotopisms are triples of the forms $(A,I,A)$ and $(I,B,B)$. All
weak associative identities such as the Bol, Moufang and extra
identities have been found to be isotopic invariant in loops for any
triple of the form $(A,B,C)$ while the central identities have been
found to be isotopic invariant only under triples of the forms
$(A,B,A)$ and $(A,B,B)$. Since associativity obeys all the
Bol-Moufang identities, the observation in the theorem agrees with
the latter stated facts.
\end{myrem}

\begin{mycor}\label{1:2}
$(L,\theta )$ is a Smarandache quasigroup with associative
subquasigroup $(S,\theta )$ if and only if any of the six Khalil
conditions is true for some four parastrophes of $(S,\theta )$.
\end{mycor}
{\bf Proof}\\
Let $(L,\theta )$ be the quasigroup in consideration. By
Lemma~\ref{0:1}, $(S, \theta )$ is a group. Notice that
$R_{s_2}R_{s_3}=R_{s_2\theta s_3}\Leftrightarrow L_{s_2\theta
s_3}^*=L_{s_3}^*L_{s_2}^*$. Hence, $(S,\theta^*)$ is also a group.
In Theorem~\ref{1:1}, two of the parastrophes are isotopes of
$(S,\theta )$ while the other two are isotopes of $(S,\theta^*)$.
Since the Khalil conditions are neccessary and sufficient conditions
for a quasigroup to be an isotope of a group, then they must be
necessarily and sufficiently true in the four quasigroup
parastrophes of $(S,\theta )$.

\begin{mylem}\label{1:3}
$(L,\theta^*)$ is a Smarandache quasigroup with associative
subquasigroup $(S,\theta^* )$ if and only if any of the following
equivalent statements is true.
\begin{enumerate}
\item $(S,\theta^*)$ is isotopic to $(S,{}^{-1}\theta )$.
\item $(S,\theta )$ is isotopic to $(S,({}^{-1}\theta )^*)$.
\item $(S,\theta^*)$ is isotopic to $(S,\theta^{-1})$.
\item $(S,\theta )$ is isotopic to $(S,(\theta^{-1})^*)$.
\end{enumerate}
\end{mylem}
{\bf Proof}\\
Replace $(L,\theta )$ with $(L,\theta^* )$ in Theorem~\ref{1:1}.

\begin{mycor}\label{1:4}
$(L,\theta^* )$ is a Smarandache quasigroup with associative
subquasigroup $(S,\theta^* )$ if and only if any of the six Khalil
conditions is true for some four parastrophes of $(S,\theta )$.
\end{mycor}
{\bf Proof}\\
Replace $(L,\theta )$ with $(L,\theta^* )$ in Corollary~\ref{1:2}.

\begin{mylem}\label{1:5}
$(L,\theta^{-1})$ is a Smarandache quasigroup with associative
subquasigroup $(S,\theta^{-1})$ if and only if any of the following
equivalent statements is true.
\begin{enumerate}
\item $(S,\theta^{-1})$ is isotopic to $(S,\theta^*)$ .
\item $(S,(\theta^{-1})^*)$ is isotopic to $(S,\theta )$.
\item $(S,\theta^{-1})$ is isotopic to $(S,{}^{-1}\theta )$.
\item $(S,(\theta^{-1})^*)$ is isotopic to $(S,({}^{-1}\theta )^*)$.
\end{enumerate}
\end{mylem}
{\bf Proof}\\
Replace $(L,\theta )$ with $(L,\theta^{-1})$ in Theorem~\ref{1:1}.

\begin{mycor}\label{1:6}
$(L,\theta^{-1} )$ is a Smarandache quasigroup with associative
subquasigroup $(S,\theta^{-1})$ if and only if any of the six Khalil
conditions is true for some four parastrophes of $(S,\theta )$.
\end{mycor}
{\bf Proof}\\
Replace $(L,\theta )$ with $(L,\theta^{-1})$ in Corollary~\ref{1:2}.

\begin{mylem}\label{1:7}
$(L,{}^{-1}\theta )$ is a Smarandache quasigroup with associative
subquasigroup $(S,{}^{-1}\theta )$ if and only if any of the
following equivalent statements is true.
\begin{enumerate}
\item $(S,{}^{-1}\theta )$ is isotopic to $(S,\theta^{-1})$.
\item $(S,({}^{-1}\theta )^*)$ is isotopic to $(S,(\theta^{-1})^*)$.
\item $(S,{}^{-1}\theta )$ is isotopic to $(S,\theta^*)$.
\item $(S,({}^{-1}\theta )^*)$ is isotopic to $(S,\theta )$.
\end{enumerate}
\end{mylem}
{\bf Proof}\\
Replace $(L,\theta )$ with $(L,{}^{-1}\theta )$ in
Theorem~\ref{1:1}.

\begin{mycor}\label{1:8}
$(L,{}^{-1}\theta )$ is a Smarandache quasigroup with associative
subquasigroup $(S,{}^{-1}\theta )$ if and only if any of the six
Khalil conditions is true for some four parastrophes of $(S,\theta
)$.
\end{mycor}
{\bf Proof}\\
Replace $(L,\theta )$ with $(L,{}^{-1}\theta )$ in
Corollary~\ref{1:2}.

\begin{mylem}\label{1:9}
$(L,(\theta^{-1})^*)$ is a Smarandache quasigroup with associative
subquasigroup $(S,(\theta^{-1})^*)$ if and only if any of the
following equivalent statements is true.
\begin{enumerate}
\item $(S,(\theta^{-1})^*)$ is isotopic to $(S,({}^{-1}\theta )^*)$ .
\item $(S,\theta^{-1})$ is isotopic to $(S,{}^{-1}\theta )$.
\item $(S,(\theta^{-1})^*)$ is isotopic to $(S,\theta )$.
\item $(S,\theta^{-1}))$ is isotopic to $(S,\theta^*)$.
\end{enumerate}
\end{mylem}
{\bf Proof}\\
Replace $(L,\theta )$ with $(L,(\theta^{-1})^*)$ in
Theorem~\ref{1:1}.

\begin{mycor}\label{1:10}
$(L,(\theta^{-1})^*)$ is a Smarandache quasigroup with associative
subquasigroup $(S,(\theta^{-1})^*)$ if and only if any of the six
Khalil conditions is true for some four parastrophes of $(S,\theta
)$.
\end{mycor}
{\bf Proof}\\
Replace $(L,\theta )$ with $(L,(\theta^{-1})^* )$ in
Corollary~\ref{1:2}.

\begin{mylem}\label{1:11}
$(L,({}^{-1}\theta )^*)$ is a Smarandache quasigroup with
associative subquasigroup $(S,({}^{-1}\theta )^*)$ if and only if
any of the following equivalent statements is true.
\begin{enumerate}
\item $(S,({}^{-1}\theta )^* )$ is isotopic to $(S,\theta )$.
\item $(S,{}^{-1}\theta )$ is isotopic to $(S,\theta^*)$.
\item $(S,({}^{-1}\theta )^*)$ is isotopic to $(S,(\theta^{-1})^*)$.
\item $(S,{}^{-1}\theta )$ is isotopic to $(S,\theta^{-1} )$.
\end{enumerate}
\end{mylem}
{\bf Proof}\\
Replace $(L,\theta )$ with $(L,({}^{-1}\theta )^* )$ in
Theorem~\ref{1:1}.

\begin{mycor}\label{1:12}
$(L,({}^{-1}\theta )^*)$ is a Smarandache quasigroup with
associative subquasigroup $(S,({}^{-1}\theta )^* )$ if and only if
any of the six Khalil conditions is true for some four parastrophes
of $(S,\theta )$.
\end{mycor}
{\bf Proof}\\
Replace $(L,\theta )$ with $(L,({}^{-1}\theta )^* )$ in
Corollary~\ref{1:2}.

\begin{myth}\label{1:13}
$(L,\pi_i )$ is a Smarandache quasigroup with associative
subquasigroup $(S,\pi_i )~\forall~i\in \{1,2,3,4,5,6\}$ if and only
if for any of some four $j\in \{1,2,3,4,5,6\}$, $(S,\pi_j)$ is an
isotope of $(S,\pi_i)$ or $(S,\pi_k)$ for one $k\in \{1,2,3,4,5,6\}$
such that $i\ne j\ne k$.
\end{myth}
{\bf Proof}\\
This is simply the summary of Theorem~\ref{1:1}, Lemma~\ref{1:3},
Lemma~\ref{1:5}, Lemma~\ref{1:7}, Lemma~\ref{1:9} and
Lemma~\ref{1:11}.

\begin{mycor}\label{1:14}
$(L,\pi_i )$ is a Smarandache quasigroup with associative
subquasigroup $(S,\pi_i )~\forall~i\in \{1,2,3,4,5,6\}$ if and only
if any of the six Khalil conditions is true for any of some four of
$(S,\pi_i )$.
\end{mycor}
{\bf Proof}\\
This can be deduced from Theorem~\ref{1:13} and the Khalil
conditions or by combining Corollary~\ref{1:2}, Corollary~\ref{1:4},
Corollary~\ref{1:6}, Corollary~\ref{1:8}, Corollary~\ref{1:10} and
Corollary~\ref{1:12}.


\begin{thebibliography}{99}
\bibitem{phd71} R. Artzy (1963), {\it Isotopy and Parastrophy of
Quasigroups}, Proc. Amer. Math. Soc. 14, 3, 429--431.
\bibitem{phd77} V. D. Belousov (1965), {\it Systems of quasigroups
with generalised identities}, Usp. Mat. Nauk. 20, 1(121), 75--146.
\bibitem{phd41} R. H. Bruck (1966), {\it A survey of binary systems}, Springer-Verlag, Berlin-G\"ottingen-Heidelberg, 185pp.
\bibitem{phd39} O. Chein, H. O. Pflugfelder and J. D. H. Smith (1990), {\it Quasigroups and Loops : Theory and Applications}, Heldermann Verlag, 568pp.
\bibitem{phd49} J. Dene and A. D. Keedwell (1974), {\it Latin squares and their applications}, the English University press Lts, 549pp.
\bibitem{phd72} R. A. Fisher and F. Yates (1934), {\it The $6\times
6$ Latin squares}, Proc. Camb. Philos. Soc. 30, 429--507.
\bibitem{phd42} E. G. Goodaire, E. Jespers and C. P. Milies (1996), {\it Alternative  Loop Rings}, NHMS(184), Elsevier, 387pp.
\bibitem{phd2} K. Kunen (1996), {\it Quasigroups, Loops and Associative Laws}, J. Alg.185, 194--204.
\bibitem{muk1} A. S. Muktibodh (2005), {\it Smarandache quasigroups rings},
Scientia Magna Journal, 1, 2, 139--144.
\bibitem{muk2} A. S. Muktibodh (2006), {\it Smarandache Quasigroups},
Scientia Magna Journal, 2, 1, 13--19.
\bibitem{phd3} H. O. Pflugfelder (1990), {\it Quasigroups and Loops : Introduction}, Sigma series in Pure Math. 7, Heldermann Verlag, Berlin, 147pp.
\bibitem{phd74} A. Sade (1959), {\it Quasigroupes parastrophiques},
Math. Nachr. 20, 73--106.
\bibitem{phd87} K. Shahbazpour (2005), {\it On identities of isotopy
closure of variety of groups}, Milehigh conference on loops,
quasigroups and non-associative systems, University of Denver,
Denver, Colorado.
\bibitem{phd76} S. K. Stein (1956), {\it Foundation of
quasigroups}, Proc. Nat. Acad. Sci. 42, 545--545.
\bibitem{phd73} S. K. Stein (1957), {\it On the foundation of
quasigroups}, Trans. Amer. Math. Soc. 85, 228--256.
\bibitem{phd75} W. B. Vasantha Kandasamy (2002), {\it Smarandache
Loops}, Department of Mathematics, Indian Institute of Technology,
Madras, India, 128pp.
\bibitem{phd83} W. B. Vasantha Kandasamy (2002), {\it Smarandache
Loops}, Smarandache notions journal, 13, 252--258.
\bibitem{phd86} W. B. Vasantha Kandasamy, {\it Smarandache
groupoids}, American Research Press (to appear).
\end{thebibliography}
\end{document}